\documentclass{amsart}
  \usepackage{amscd,amssymb,epsfig}
   \usepackage{epic,eepic}

                     \numberwithin{equation}{subsection}

                     \newtheorem{propo}{Proposition}[subsection]
                     \newtheorem{corol}[propo]{Corollary}
                     \newtheorem{theor}[propo]{Theorem}
                     \newtheorem{lemma}[propo]{Lemma}
                     \theoremstyle{definition}

                     \theoremstyle{remark}

		     \newcommand{\E}{\mathcal{A}}\newcommand{\MM}{\mathcal{Q}}
		     
		     \newcommand{\QQ}{\mathbb{Q}}
                     \newcommand{\ZZ}{\mathbb{Z}}
                     \newcommand{\RR}{\mathbb{R}}

		     \newcommand{\A}{\mathcal{A}}

                     \newcommand{\Hom}{\operatorname{Hom}}
		     \newcommand{\Ker}{\operatorname{Ker}}
		     \newcommand{\Exp}{\operatorname{Exp}}

\newcommand{\T}{\mathcal{T}}
\newcommand{\V}{\,\vert\,}
\newcommand{\PP}{\mathcal{P}}\newcommand{\WW}{\mathcal{W}}
\newcommand{\VV}{\mathcal{V}}
\newcommand{\UU}{\mathcal{U}}

		      \newcommand{\card}{\operatorname{card}}
		      
                     \newcommand{\id}{\operatorname{id}}

		    \newcommand{\modu}{\operatorname{mod}}

		  \newcommand{\Perm}{\operatorname{Perm}}\newcommand{\End}{\operatorname{End}}

\begin{document}
      \title{Coalgebras of   words and    phrases}
                     \author[Vladimir Turaev]{Vladimir Turaev}
                     \address{%
              IRMA, Universit\'e Louis  Pasteur - C.N.R.S., \newline
\indent  7 rue Ren\'e Descartes \newline
                     \indent F-67084 Strasbourg \newline
                     \indent France \newline
\indent e-mail: turaev@math.u-strasbg.fr }
                     \begin{abstract} We  introduce two  constructions of a  coassociative   comultiplication  in the algebra    of phrases in a given alphabet. As a preliminary step	we give two  constructions of a  pre-Lie comultiplication in the module   generated by words.   
                     \end{abstract}
                     \maketitle

  \section {Introduction}

We give two new   constructions of  non-commutative non-cocommutative  Hopf algebras  (of infinite rank). These constructions  can be summarized as follows. A module $M$ over a commutative 
ring gives rise to the tensor algebra $\WW=T(M)=\oplus_{m\geq 0} M^{\otimes m}$ and to its positive subalgebra $\VV=T_+(M)= \oplus_{m\geq 1} M^{\otimes m}$. Under  certain additional assumptions, we  define       comultiplications in the tensor  algebras $\PP=T(\VV)$ and   $  \MM =T(\WW) $ so that  they  become  Hopf algebras.  

 A comultiplication $\Delta$ in a  graded algebra $\oplus_{m\geq 0} T^m$    has a  leading term  which is the homomorphism $T^1\to T^1\otimes T^1$ obtained by applying $\Delta$ and then projecting to $T^1\otimes T^1$. 
This leading term  yields the first approximation to $\Delta$ and  is often interesting  in itself.   The leading terms of our  comultiplications in $\PP $ and   $  \MM   $ are homomorphisms $\VV\to \VV \otimes \VV$ and $\WW\to \WW \otimes \WW$. They turn out to be   pre-Lie comultiplications in the sense explained below. In particular, dualizing and skew-symmetrizing them, we obtain Lie brackets in the dual modules $\VV^*$ and  $\WW^*$.  
  
The study of tensor algebras   can be   reformulated in terms of words and phrases. Suppose from now on  that the module $M$ is  free   with  basis $\E$. Consider the  basis of   $\VV=T_+(M)$ formed by the vectors  $A_1\otimes A_2\otimes...\otimes A_m$ where $m\geq 1$ and $A_1, A_2,...,A_m\in \E$. Omitting the symbol $\otimes$ we can identify these vectors with  words in the alphabet $\E$.  A basis in $\PP= T(\VV)$ is formed by finite sequences of basis vectors in $\VV$, that is by sequences of words. We call such sequences {\it phrases}.  Natural bases in $\WW=T(M)$ and $\MM=T(\WW)$ are similarly interpreted  as words and phrases with the difference that here we allow an empty word.  This  interpretation    gives a pleasant linguistical flavour to the theory.  It   places the study of   words and  phrases   in the  setting 
  of   Hopf algebras, Lie algebras,  and related  algebraic objects. 
  
One well known  construction of a comultiplication  in   $  \WW$   is based on the notion of a subword. For us, a {\it word}   is a  finite  sequence  of elements of $\E$ and a {\it subword} is any subsequence. Given a word $w$ we can form the sum $\sum   w'\otimes w/w'$ where
$w'$ runs over all subwords of $w$ and the word $w/w'$ is obtained by deleting $w'$ from $w$. This   yields  the familiar {\it shuffle   comultiplication} in $\WW$. It is coassociative and cocommutative. 

Our   comultiplications   are  based on two  different  constructions. 
They  are determined by certain additional data which should be fixed from the very beginning. The comultiplications in $\VV$ and $\PP=T(\VV)$
 depend  on a choice of a so-called stable set of words $L$ in the alphabet $\E$.   We associate with a word $w$ the  sum $\sum   w'\otimes w/w'$ where $w'$ runs over all   subwords of $w$ formed by   consecutive  letters  and belonging to $L$.  This  yields  a pre-Lie comultiplication $\rho_L$ in $\VV$. In a similar way we construct a coassociative  comultiplication    in $\PP$ with leading term $\rho_L$. 
  
  The second construction  begins with   fixing a   mapping $\mu$ from $\E\times \E$ to the ground ring. The corresponding pre-Lie comultiplication $\rho_\mu$ in $\WW$ is defined as follows. For a  
subword  of length two $a=AB$ of  a    word $w$, set $\mu_a=\mu(A,B)$ and denote by  $w'_a$ be the subword of $w$ formed by   the letters    of $w$ appearing   between $A$ and $B$. Let $w''_a$ be the word obtained by deleting both $a$ and $w'_a$ from $w$. Then $ \rho_\mu (w)= \sum_a\mu_a w'_{a}\otimes w''_a$. In a similar way we construct   a coassociative  comultiplication   in $\MM=T(\WW)$ with leading term $\rho_\mu$. 
 
Under a certain   choice of $\mu$, the  latter  comultiplication   is closely  related to the Connes-Kreimer    comultiplication in the algebra of rooted trees, see \cite{ck}, or, more precisely, to its non-commutative version for planar rooted trees due to Foissy \cite{fo}.

  Although we can directly define our   comultiplications in $\PP$ and $\MM$, we   begin with a study of their leading terms.  In Section 2 we discuss   relevant notions from the theory of pre-Lie multiplications and comultiplications. Sections 3 and 4 are concerned with  the   comultiplication in   $\PP$   derived from a stable set  of words: the leading term is defined in  Section 3 and  the  comultiplication itself is defined in Section 4.
Sections 5 and 6 are concerned with  the comultiplication   in  $\MM$ derived from a mapping $\mu$ as above: the leading term is defined in  Section 5 and  the   comultiplication   itself is defined in Section 6. At the end of Section 6 we discuss connections with the theory of planar rooted trees.

 Throughout the paper, we fix a
commutative ring
 with unit $R$. The symbol $\otimes$ denotes the tensor product of 
$R$-modules over $R$. 

 \section{Pre-Lie  coalgebras and left-handed bialgebras}\label{cob}

   \subsection{Pre-Lie algebras}\label{su1}  By a {\it multiplication} in an $R$-module $M$ we   mean an $R$-bilinear mapping  $M\times M \to M$.
A (left) {\it  pre-Lie  algebra}  over $R$  is  an $R$-module $M$ endowed 
with a  multiplication  $M\times M \to M$, denoted
$\star$, such that for any
$f,g,h\in M$,
	 \begin{equation}\label{jacla}   f \star 
(g \star  h) -(f\star  g) \star  h  =   g \star  (f \star 
h) -(g\star  f) \star  h .\end{equation}
The mapping $\star$ is called then a {\it pre-Lie  multiplication}, see   \cite{ge} and \cite{vi}.
For example,  any associative multiplication is pre-Lie. A fundamental property of a pre-Lie multiplication 
$\star$ in $M$ is that the formula $[f,g]=f\star g-g \star f$ defines a Lie bracket in $M$. The Jacobi identity
$$ [f, [g,h]]+[g,[h,f]]+ [h,[f,g]]=0 $$
 for    $f,g,h\in M$ is a direct consequence of (\ref{jacla}). In this way every pre-Lie algebra becomes a Lie algebra.

There is a similar notion of right pre-Lie algebras. We will consider 
only left pre-Lie algebras and refer to them simply as pre-Lie 
algebras. Non-trivial examples of pre-Lie 
algebras can be obtained from derivations in   algebras.   A  {\it  derivation} in an associative $R$-algebra $M$ is an  $R$-linear homomorphism $d:M\to M$ such that 
 $d(fg)= d(f) g+ f \,d(g)$ for any   $f,g\in M$. An easy computation shows that for any $a,b\in M$, the formula $f\star g=a fg +bf  d(g) $ defines a pre-Lie  multiplication in $M$.

\subsection{Pre-Lie   coalgebras}\label{foplsu1}   A  {\it comultiplication} in an $R$-module $A$ is an $R$-linear  homomorphism $\rho: 
A\to    A \otimes A$. We   associate with  such $\rho$ a homomorphism $\tilde \rho: A\to    A \otimes A\otimes A$ sending $a\in A$ to  $$\tilde \rho (a)
= (\id_A 
\otimes \rho) \rho (a)-  (\rho \otimes \id_A) \rho (a) .$$ The comultiplication $\rho$ is {\it coassociative} if $\tilde \rho =0$.

Dualizing   Formula  \ref{jacla}, we obtain a notion of    a 
pre-Lie coalgebra. Namely, a   {\it pre-Lie coalgebra}
is an $R$-module $A$ endowed with a  comultiplication $\rho: 
A\to A^{\otimes 2}=A\otimes A $ such that for all $a\in A$,
 \begin{equation}\label{prejacla}  P^{1,2}  \,\tilde \rho (a)= \tilde \rho (a).  \end{equation}
   Here and below given a  module $U$, we denote by $P^{1,2} $   the
 endomorphism of $ {A^{\otimes 2}}\otimes U$ permuting the  first two 
tensor factors, i.e., mapping $a\otimes b \otimes u$ to $b\otimes a \otimes u $ for $a,b \in A, u\in U$.     A  comultiplication $\rho$ satisfying (\ref{prejacla}) is   a {\it pre-Lie comultiplication}.
Clearly, a  coassociative comultiplication is pre-Lie.

Given   a pre-Lie coalgebra $(A, \rho)$, consider the   module 
$A^*=\Hom_R(A,R)$ and  the   evaluation pairing $\langle \,, \,\rangle: A\otimes A^* \to R$. The comultiplication $\rho$   induces a   pre-Lie multiplication $\star_\rho$ in $A^*$   by  
 \begin{equation}\label{jcccc}  \langle a,  f\star_\rho  g\rangle=\sum_i \langle a'_i, f\rangle \,  \langle a''_i, g\rangle  \in R\end{equation}
for any  $f,g \in A^*$, $a\in A$  and any (finite) expansion $\rho (a)=\sum_i a'_i \otimes 
a''_i\in A^{\otimes 2}$. In this way $A^*$ becomes a pre-Lie algebra.

\subsection{Comodules over pre-Lie coalgebras}\label{su3}  A {\it (left) comodule} over a pre-Lie coalgebra $(A,\rho)$
is a pair (an $R$-module $U$, an $R$-homomorphism $\theta: U\to A\otimes U$) 
such that   $\tilde \theta= (\id_A  \otimes \theta) \theta - (\rho \otimes \id_U)\theta
:U\to A^{\otimes 2}\otimes U$ satisfies  
$ P^{1,2}  \tilde \theta  = \tilde \theta$. 
 Given a comodule $(U,\theta)$ over   $(A,\rho)$, we  define a right action of $A^*$ on $U$ by
   $$uf=   \sum_{i} \langle u'_i, f\rangle \, u''_i$$  
   for any  $f  \in A^*$, $u\in U$,  and any   expansion $\theta (u)=\sum_i u'_i \otimes 
u''_i $ with $u'_i \in A, u''_i \in  U$.

\begin{lemma}\label{colbra}  The action of $A^*$ on $U$ is a   Lie algebra action.
\end{lemma}
                     \begin{proof} We must prove that  $u[f,g] = (uf)g -(ug)f$ for all  $f,g \in A^*$, $u\in U$ where $[f,g]=f\star_\rho g-g \star_\rho f$.  We      shall use Sweedler's notation  for the expansion    $\theta (u)=\sum_i u'_i \otimes u''_i$
    and  write it in the form $\theta (u)=\sum_{(u)} u' \otimes u''$. Similarly, for $a\in A$, we write $\rho (a)=\sum_{(a)} a' \otimes a''$. 
 By definition,  $$(uf)g=   \sum_{(u)} \langle u', f\rangle u''g  
		     =  \sum_{(u), (u'')} \langle u', f\rangle \langle (u'')', g\rangle (u'')'' .$$
		     Hence
		     $$ (uf)g -(ug)f =\sum_{(u), (u'')}  \langle u', f\rangle \langle (u'')', g\rangle (u'')''-\langle (u'')', f\rangle \langle u', g\rangle (u'')''.$$
		     To compute $  u[f,g]    $, observe  that for $a\in A$, 
		     $$\langle a, [f,g]\rangle=\langle a, f \star_\rho g -g \star_\rho f\rangle
		      =\sum_{(a)} \langle a', f\rangle \,  \langle a'', g\rangle - \langle a'', f\rangle\, \langle a', g\rangle .$$
		     Therefore 
		     $$  u [f,g]   =  \sum_{(u)} \langle u', [f,g]\rangle\,  u''
		     =  \sum_{(u), (u')}  \langle (u')', f\rangle \,   \langle (u')'', g\rangle\,  u'' - \langle (u')'', f\rangle \, \langle (u')', g\rangle \, u''
		.$$
		     The  equality $u[f,g] = (uf)g -(ug)f$   follows now 
		     from the formula
		     $$\sum_{(u), (u')}  (u')'\otimes (u')''\otimes u''
		     -(u')''\otimes (u')'\otimes u''=\sum_{(u), (u'')} u'\otimes (u'')'\otimes (u'')'' - (u'')'\otimes u'\otimes (u'')''$$
		     which  is a reformulation of the equality  
$ P^{1,2}  \tilde \theta  = \tilde \theta$. 
		\end{proof}  

Lemma \ref{colbra} implies that the formula $fu=-uf$ with $u\in U, f\in A^*$ defines a left Lie algebra action of $A^*$ on $U$. 

It is clear that   $(A, \theta=\rho:A\to A\otimes A)$ is a comodule over $(A,\rho)$.      Lemma  \ref{colbra} yields a   Lie algebra action of $A^*$  on $A$.  

\subsection{Left-handed bialgebras}\label{su4} By a {\it bialgebra} we shall mean a pair $(T,\Delta)$ where $T$ is an associative unital $R$-algebra and $\Delta:T\to T\otimes T$ is a coassociative  algebra comultiplication   in   $T$.
The words \lq\lq algebra comultiplication" mean  that  $\Delta(1)=1\otimes 1$ and
$  \Delta(ab)=   \Delta(a)     \Delta(b)$ for any $a,b\in T$. Here 
  multiplication in $T \otimes T$ is defined  by   $(a\otimes a') (b \otimes b')= ab\otimes a'b'$ for $a,b,a',b'\in T$.  We do not require the existence of a counit $T\to R$ although in all our constructions of bialgebras there will be a counit.
  
The  comultiplication   in a bialgebra $(T, \Delta)$ induces an associative multiplication     in    $T^*=\Hom_R (T,R)$
by $fg (a) =\sum_{(a)} f(a')\, g(a'')$ for $a\in T, f,g\in T^*$, and  any   expansion $\Delta (a)=\sum_{(a)}  a' \otimes a''$. This makes $T^*$ into an associative algebra. If $T$ has a counit, then $T^*$ is a unital algebra.
 Dualizing      multiplication in $T$ we obtain a    homomorphism  from $T^*$
to a certain completion of $T^*\otimes T^*$.  We  call this homomorphism  {\it quasi-comultiplication}.

A bialgebra $(T,\Delta)$  is {\it graded} if  $T$ splits as a direct sum of submodules  $T=   \oplus_{m\geq 0} T^m$ such  that $1\in T^0$ and $T^m T^n\subset T^{m+n}$ for all $m,n$. Set $T_+= \oplus_{m\geq 1 } T^m \subset T$. 
 A graded bialgebra $(T,\Delta)$
 is   {\it left-handed} if for any $a\in T^1 $,  
   $$\Delta (a)-a\otimes 1 -1\otimes a\in   T_+\otimes T^1 .$$  The {\it leading term} of   $\Delta$   is then  the homomorphism  
   $(\pi \otimes \pi ) \Delta\vert_{T^1} : T^1\to    T^1\otimes T^1$ where 
    $\pi:T\to T^1$ is the   projection.    
   We shall mainly apply these definitions in the case where 
  $$T=T(A)=\oplus_{m\geq 0} A^{\otimes m}$$ is the tensor algebra
  of an $R$-module $A$. 
Here    $A^{\otimes 0}=R$, $A^{\otimes 1}=A$, 
and $A^{\otimes m}$
is the
  tensor product of $m$ copies  of $A$ for $m\geq 2$.

The next lemma relates left-handed bialgebras with pre-Lie coalgebras.

	 \begin{lemma}\label{12covercoc} Let $(T=   \oplus_{m\geq 0} T^m, \Delta)$ be a  left-handed graded  bialgebra  such that   $T^1$ generates $T_+= \oplus_{m\geq 1} T^m$.   Then the leading term $\rho  : T^1\to    T^1\otimes T^1$    of $\Delta$ is a pre-Lie  comultiplication in $T^1$.   
  \end{lemma}

  \begin{proof}    For $a\in T^1$, we can expand
		  $$\Delta(a)= a\otimes 1+1\otimes a+ \sum_i a_i'\otimes a_i''+\sum_j a_j^1a_j^2 \otimes a_j^3 \,\, (\modu \oplus_{m\geq 3 } T^m\otimes T^1)$$ 
	where $i,j$ run over finite sets of indices and 
		  $a_i', a_i'',a_j^1,a_j^2,a_j^3 \in T^1$.     
		  Clearly,
		  $\rho (a) =\sum_i a'_i\otimes a''_i$.  Our assumptions imply that $\Delta (T^m)\subset \oplus_{k+n\geq m} T^k\otimes T^n$ for all $m$.
		  Computing  $(\id \otimes \Delta) \Delta(a) $ and
		  $(\Delta\otimes \id) \Delta(a) $
  modulo $\oplus_{k,n,r, k+n+r\geq 4} \,T^k \otimes T^n\otimes T^r$  and equating the resulting expressions
  we obtain 
  that $$\tilde \rho (a)=   \sum_j a^1_j\otimes a^2_j
\otimes a^3_j  +   a^2_j\otimes a^1_j \otimes a^3_j. $$
 Therefore  $  P^{1,2}    \tilde \rho (a)=\tilde \rho (a)$.
  \end{proof}

   Let $(T , \Delta)$ be as in the lemma. A {\it (left) comodule} over   $(T, \Delta)$ is a pair (an $R$-module $U$, an $R$-homomorphism $\Theta:U\to T\otimes U$) such that $(\Delta\otimes \id_U) \Theta
 =(\id_T\otimes \Theta) \Theta$ and $\Theta (u)- 1\otimes u\in   T_+\otimes U $ for any $u\in U $.   Then $U$ becomes a right  module over the   algebra $T^*$
 by  $uf=\sum_{(u)} f(u') u''$ for $f\in T^*$, $u\in U$, and any finite expansion $\Theta(u)= \sum_{(u)}  u' \otimes u''$.   For example, set $U=T^1$ and define $\Theta: U\to T\otimes U$   by   $ \Theta(u)=\Delta(u)- u\otimes 1  $ for $u\in U$. The pair $(U,\Theta)$ is a comodule over  $(T, \Delta)$ (cf. the proof of  Theorem \ref{comumu} below). This makes $U=T^1$ into      a     right $T^*$-module.

   \subsection{Remarks}\label{exxx} 1. Lemma \ref{12covercoc} suggests to study    when a given pre-Lie comultiplication is induced from a left-handed  comultiplication in a graded bialgebra. For the pre-Lie comultiplications defined in the next sections this will be always the case.

    2. Consider a pre-Lie coalgebra $(A, \rho)$ and  a submodule   $B\subset A$ such that
 $\rho(A)\subset B\otimes A$ and $\rho(B)\subset B\otimes B$. Then $(B, \rho\vert_{B})$ is a pre-Lie coalgebra  and      $(A,\theta=\rho:A\to B\otimes A)$ is a comodule over   $(B, \rho\vert_{B})$. The associated action  of $B^*$ on $A$ is compatible with the action of $A^*$  on $A$  via the Lie algebra homomorphism $A^*\to B^*$ induced by the inclusion $B\subset A$.  

3. For a   pre-Lie coalgebra $(A, \rho)$, the homomorphism  $  \rho- P^{1,2} \rho:A\to A^{\otimes 2}$ is a Lie 
cobracket (cf., for instance \cite{tu1}). In particular in Lemma \ref{12covercoc}, $\rho-P^{1,2}\rho$ is a Lie cobracket in $T^1 $. This still holds   if the left-handedness assumption on $\Delta$ is weakened to    $\Delta (a)-a\otimes 1 -1\otimes a\in T_+ \otimes T_+$ for all $a\in T^1$. 
      
4. Lemma \ref{12covercoc} extends to comodules  as follows.  Let $(T, \Delta)$ be as in this lemma.  The leading term of  a comodule $(U,\Theta)$   over   $(T, \Delta)$  is   the homomorphism  
   $\theta=(\pi \otimes \id ) \Theta  : U\to    T^1\otimes U$ where 
    $\pi:T\to T^1$ is the   projection.      A direct computation shows that    $(U,\theta)$ is a comodule over the pre-Lie coalgebra $(T^1,\rho=(\pi \otimes \pi ) \Delta\vert_{T^1})$. Note also that if  $T^0=R$, then   the projection $T\to T^0=R$ is a counit of $(T,\Delta)$.

   \section{Coalgebra   of words}\label{coac}

 \subsection{Words}\label{modic} For us, an {\it alphabet} is an arbitrary set
  and {\it letters } are its elements. Throughout the paper we fix an alphabet $\E$.  A  {\it word of length $m\geq 1$} is  a mapping from the set $\{1,2,...,m\}$ to    $\E$. By definition, there is a unique word of length $0$ called the 
 {\it empty word} and denoted $\phi$.  To  present a word $w:\{1,2,...,m\}\to \E$ it is enough to  write down the sequence of letters $w(1) w(2)...w(m)$. For instance, the symbol $ABA$ represents the word $\{1,2,3\}\to \E$ sending $1$ and $3$ to $A\in \E$ and sending $2$ to $B\in \E$. 
 
 Writing down consecutively the letters of two  words $w$ and   $x$ we obtain their concatenation $wx$. For instance,   the concatenation of  $w=ABA$ and $x=BB$ is the word $wx=ABABB$.
 
 A word $x$ is a {\it factor} of a word $w$ if $w=yxz$ for certain words $y,z$. A factor $x$ of $w$ is {\it proper} if $x\neq w$.
 
Given  a word $w $ of length $m\geq 1$ and   numbers $1\leq i\leq j\leq m+1$, set $w_{i,j}=w(i) w(i+1)...w(j-1)$. This   word of length $j-i$  is a factor  of $w =w_{1,i}w_{i,j}w_{j, m+1}$. The word $w_{i,j}$ is empty iff $i=j$ 
    and proper iff $(i,j)\neq (1,m+1)$.

   Let $\WW=\WW(\E)$ be the free $R$-module freely  generated by the set of   words in the alphabet  $\E$. A  typical element of $\WW$ is a finite formal linear combination of   words with coefficients in $R$. Each  word $w$ 
 represents a vector in $\WW$ denoted also $w$. These vectors     form a basis of $\WW$. 
Let $\VV=\VV(\E)$ be the submodule   of $\WW$   generated by non-empty words. Clearly, $\WW=\VV\oplus R\phi$. 
 
\subsection{Stable sets of words}\label{lumodic}  A set $L$ of non-empty words  is    {\it stable} if it satisfies   the following    condition:
  
$(\ast)$ For any word $w$   of length $m\geq 1$ and any   indices $1\leq i< j \leq m+1$   such that $(i,j)\neq (1,m+1)$
and   $w_{i,j}\in L$,  the word $w$ belongs to $L$ if and only if the word $w_{1,i}   w_{j,m+1}$ belongs to $L$.

The words belonging to a stable set $L$ will be called {\it $L$-words}.
The \lq\lq only if" part of $(\ast)$  means that   striking out from any $L$-word  a   proper factor belonging to $L$  one obtains   an $L$-word.   The \lq\lq  if" part of $(\ast)$ means that   inserting an    $L$-word in   an  $L$-word one obtains  an $L$-word. In particular,  concatenation of   $L$-words is an $L$-word.
 
We give a few examples of stable sets of words.

(1) The set of all non-empty words  and the void set of words are stable.

To give the next example, denote   by $f_A(w)$ the number of appearances of a letter $A\in \E$ in a word $w$. For instance,  $f_A(ABA)=2$. 

(2) Pick   $A\in \E$. The  set of non-empty words $w$ such that $f_A (w)=0$ is stable. Similarly, for any integer $N$, the set of non-empty words $w$ such that $f_A (w)$ is divisible by $N$ is stable.

(3) Pick   letters $A_1,..., A_n\in \E$ and   elements $g_1,..., g_n $ of a certain abelian group.
The  set of non-empty words $w$ such that $  f_{A_1} (w) g_1+...+   f_{A_n} (w) g_n=0$ is stable. 

Observe finally that the intersection of any family of stable sets of words is stable. A union of stable sets can  be non-stable.

\subsection{Pre-Lie coalgebra of   words}\label{prelilumodic}  
 We fix a  stable set of words $L$ and derive from it  a pre-Lie comultiplication $\rho_L$  in the module $\VV$.   
 
A {\it simple  cut} of a word $w$ of length $m\geq 1$ is a pair of indices $1\leq i< j\leq m+1$   such that $(i,j) \neq (1,m+1)$ and $w_{i,j}\in L$. To indicate that $(i,j)$ is a simple cut of $w$ we write $(i,j)\prec w$. Set
 $$\rho_L(w)=\sum_{(i,j)\prec w} w_{i,j} \otimes  w_{1,i} w_{j,m+1}\in \VV\otimes \VV $$
 where the sum runs over all simple  cuts of $w$. 
 Note that the words $w_{i,j}$ and $ w_{1,i} w_{j,m+1}$ are necessarily non-empty. 
   This defines $\rho_L$ on the basis of $\VV$ and extends by linearity to a comultiplication $\rho_L:\VV\to \VV\otimes \VV$. 
For example, if $A,B\in \E$ and $L$ is the set of all non-empty words, then $\rho_L(A)=0$, $\rho_L (AB)=A\otimes B +B\otimes A$, 
$$\rho_L (ABA)=  A\otimes BA   +   B\otimes AA  +A \otimes AB+ AB \otimes A+BA \otimes A .$$

  \begin{theor}\label{epr}  $\rho_L$    is a pre-Lie comultiplication in $\VV$.   
  \end{theor} 
  \begin{proof} For a word $w$ of length $m\geq 1$, 
$$(\id\otimes \rho_L) \rho_L (w)=  \sum_{(i,j)\prec w} w_{i,j} \otimes \rho_L( w_{1,i} w_{j,m+1}).$$
To compute $w_{i,j} \otimes \rho_L(w_{1,i} w_{j,m+1})$, we need to consider  factors
of $w_{1,i} w_{j,m+1}$. There are three kinds of them: (i) the factors of $w_{1,i-1}$, (ii)  the factors of $w_{j+1,m+1}$, and (iii) the factors obtained by concatenation of some $w_{i',i}$ with some $w_{j,j'}$ where $i'\leq i, j\leq j'$. Here and below     $i,j,i',j'\in \{1,2,...,m\}$. The factors of $w_{1,i} w_{j,m+1}$ of type (i) contribute  to $w_{i,j} \otimes \rho_L(w_{1,i} w_{j,m+1})$ the   expression
$$x_{i,j}=  \sum_{(i',j')\prec w, j'< i} w_{i,j} \otimes w_{i',j'} \otimes w_{1,i'} w_{j',i} w_{j,m+1}.$$
The factors of $w_{1,i} w_{j,m+1}$ of   type (ii) contribute  
$$y_{i,j}=  \sum_{(i',j')\prec w, j <  i' } w_{i,j} \otimes w_{i',j'} \otimes w_{1,i} w_{j,i'} w_{j',m+1}.$$
The factors of $w_{1,i} w_{j,m+1}$ of    type (iii)  contribute  
$$z_{i,j}= \sum_{ i'\leq i, j\leq j', (i',j')\neq (1,m+1),  w_{i', i } w_{j ,j'} \in L} w_{i,j} \otimes w_{i',i} w_{j,j'} \otimes  w_{1,i'} w_{j',m+1}.$$
Recall that $w_{i,j} \in L$. By definition of a stable set of words,   $w_{i', i} w_{j,j'} \in L$ iff $(i',j')\neq (i,j)$ and $w_{i',  j'} \in L$. Therefore
the conditions on $i',j' $ in the latter sum are equivalent to the conditions $i'\leq i, j\leq j', (i',j')\prec w,  (i',j')\neq (i,j)$. We obtain
$$\sum_{(i,j)\prec w} z_{i,j}=\sum_{(i,j)\prec w, (i',j')\prec w, i'\leq i, j\leq j', (i',j')\neq (i,j)} w_{i,j} \otimes w_{i',i} w_{j,j'} \otimes  w_{1,i'} w_{j',m+1}
$$
$$=\sum_{(i',j')\prec w}\,\,\rho_L( w_{i',j'})   \otimes  w_{1,i'} w_{j',m+1}=(\rho_L\otimes \id) \rho_L (w).$$
Hence, $$\tilde \rho_L (w)= (\id\otimes \rho_L) \rho_L(w)-(\rho_L\otimes \id) \rho_L(w)= \sum_{(i,j)\prec w} x_{i,j}+ y_{i,j}.$$  It remains to observe that
$$P^{1,2}( \sum_{(i,j)\prec w} x_{i,j})= \sum_{(i,j)\prec w} y_{i,j} .$$ 
Therefore  $\tilde \rho_L (w)$ is invariant under $P^{1,2}$. 
  \end{proof}

\subsection{Word indicators}\label{fomph} The elements of the module  $\VV^*=\Hom_R(\VV,R)$ are called {\it word indicators}.  The module $\VV^*$ admits a decreasing filtration $\VV^*=   \VV^{*(1)}\supset \VV^{*(2)}\supset ...$ where $\VV^{*(m)}$  consists of the word indicators annihilating all non-empty words of length $\leq m-1$. Clearly, $\cap_m  \VV^{*(m)}=0$.
We can  consider infinite sums in  $\VV^*$ as follows. Let $f_1,f_2,\ldots
\in \VV^*$
be a  sequence of word indicators such that for any $m\geq 1$ all terms of the
sequence
starting from a certain place belong to $\VV^{*(m)}$. Then for any $a\in \VV$,
the sum $f(a)=f_1(a)+f_2(a)+\cdots$ contains only a finite number of
non-zero
terms and the formula $a\mapsto f(a):\VV\to R$   defines
a word indicator    $f=f_1+ f_2+ \cdots$.  A similar construction shows that
the natural homomorphism   $\VV^*\to \projlim_m 
(\VV^*/\VV^{*(m)})$
is an
isomorphism.

 By the general theory, 
 the pre-Lie comultiplication $ \rho_L$ in $\VV$ induces a pre-Lie multiplication $ \star_L$  and a Lie bracket $[f,g]_L=f\star_L g- g\star_L f$ in   $\VV^* $. 
 It is clear that $\VV^{*(m)} \star_L \VV^{*(n)}\subset \VV^{*(m+n)}$ and  $[\VV^{*(m)} , \VV^{*(n)}]_L\subset \VV^{*({m+n})}$ for all $m,n $. 
 Therefore $ \VV^{*} =\projlim_m 
(\VV^*/\VV^{*(m)})$ is a projective limit of nilpotent Lie algebras.  
 
 Recall the  Lie algebra action   of 
  $\VV^*$  on $\VV$ induced by $\rho_L$. For   $f\in \VV^*$ and a word $w$ of length $m\geq 1$,
 $$fw=-wf=-\sum_{(i,j)\prec w}  \langle w_{i,j}, f\rangle \,    w_{1,i} w_{j,m+1}.$$
All word indicators annihilating $L$ lie in the kernel of this action.

 The module $\VV$ has an increasing filtration $0 =\VV_0\subset \VV_1\subset \VV_2\subset ...$
 where $\VV_m$ is generated by   non-empty words of length $\leq m$. It is clear that   $f\VV_m \subset \VV_{m-1}$ for all $m $ and    $f\in \VV^*$. Thus the action of $\VV^*$ on $\VV_m$ is nilpotent for all $m$. 
 
The restriction of  $\rho_L:\VV\to \VV\otimes \VV$  to the  module $ RL\subset \VV$  generated by $L$ is a pre-Lie comultiplication in $RL$ since $\rho_L(RL)\subset RL\otimes RL$. Clearly $\rho_L(\VV)\subset RL\otimes \VV$  so that  $\VV$ becomes a  comodule over the pre-Lie coalgebra $RL$.  The associated actions  of $(RL)^*$ and $\VV^*$ on $\VV$ are compatible   via the Lie algebra homomorphism $\VV^*\to (RL)^*$ induced by the inclusion $RL\subset \VV$.  This is a special case of Remark \ref{exxx}.2.

\subsection{Example}\label{exidifr} Let the alphabet $\E$ consist  of only one letter $A$. Denote by $A^m$ the   word consisting of $m$ letters $A$ where $m\geq 1$. A stable set of words $L$ is determined by a positive integer $N$ and  consists of  all words  $A^m$ with $m$ divisible by $N$. The   comultiplication $\rho_L$ is  computed by
 $$\rho_L(A^m)=\sum_{k\geq 1, kN< m} (m+1-kN)  A^{kN} \otimes A^{m-kN} $$
 for any $m\geq 1$. 
We can  identify a word indicator $f$  with  its  generating function $ \sum_{m\geq 1} f(A^m) t^m$. One easily computes that
 $f\star_L g= (f\V_N) \cdot (g+tg')$ where for a formal power series $f=\sum_{m\geq 1} a_m t^m$, we set $f\V_N =\sum_{m\geq 1} a_{mN}  t^{mN}$ and $f'=\sum_{m\geq 1} m a_{m}  t^{m-1}$.
   For $N=1$, we obtain   $ f\star_L g = fg+tfg'$.

 \subsection{The group $\Exp_L \VV^*$}\label{succ3}  Assume in this subsection that $R\supset \QQ$. 
 For   $f,g\in \VV^*$,
set
$$fg= f+g + \frac{1}{2} [f,g]_L+ \frac{1}{12} ([f,[f,g]_L]_L +[g,
[g,f]_L]_L)+
\cdots \in \VV^*$$
where the right-hand side is the Campbell-Baker-Hausdorff series for $\log
(e^f e^g)$. The resulting mapping $  \VV^*\times \VV^*\to \VV^*$ is a
group
multiplication in $\VV^*$. Here $f^{-1}=-f$ and $0\in  \VV^*$  is the group unit.
The resulting group
  is denoted $\Exp_L \VV^*$. Heuristically, this is
the \lq\lq Lie group" with Lie algebra $(\VV^*, [\,,\,]_L)$. The equality  $\VV^*= 
\projlim_m (\VV^*/\VV^{*(m)})$
implies that $\Exp_L \VV^*$  
is   
pro-nilpotent.   

The   action of $\VV^*$ on  $\VV$ induced by $\rho_L$ integrates  to a group action of $\Exp_L \VV^*$ on $\VV$. To see this,  denote by $\varphi(f)$ the additive endomorphism $a\mapsto fa$ of $\VV$ determined by $f\in  \VV^*$.  Set 
$$e^{\varphi(f)}(a)=\sum_{k\geq 0} \frac {1} {k!} (\varphi (f))^k (a).$$
The latter sum makes sense since for any   $a\in \VV$, only a finite number of   terms   in the sum are non-zero.
  The formula  $f\mapsto e^{\varphi(f)} $ defines a group homomorphism $ \Exp_L \VV^*\to \End (\VV)$, i.e.,  a group action of 
$\Exp_L \VV^*$ on $\VV$.
 
 \subsection{Remarks}\label{ipisl} 1. There is an embedding  $\delta: \VV\hookrightarrow \VV^*$  mapping a word $w$ into the word indicator $\delta_w$ whose value on $w$ is 1 and whose value on all other words is $0$. It is easy to check  that the image of $\delta$ is closed under $\star_L$. This induces a pre-Lie  multiplication $ \circ_L$ in   $\VV$. For words $w=A_1A_2... A_m$ and $x=B_1 B_2...B_n$ with $m,n\geq 1$, we have $w\circ_L x=0$ if $w\notin L$ and 
$$w\circ_L x =\delta^{-1} (\delta_w \star_L \delta_x) =\sum_{i=0}^n B_1 ...B_i A_1 ... A_m B_{i+1}...B_n $$
if $w\in L$.
In the case where $L$ is the set of all non-empty words, this pre-Lie multiplication in $\VV$ is essentially due to Gerstenhaber \cite{ge}.

2. With a slight modification, the notion of a stable set  of words can be used to define a coassociative  comultiplication in   $\WW=\VV\oplus R \phi$. Let us say that a set of words $S$ is {\it strongly stable} 
if $\phi\in S$ and for any word $w$ and any its subword $w'$ we have $w'\in S \Rightarrow (w\in S \Leftrightarrow w/w'\in S)$ where $w/w'$ is the word obtained by deleting $w'$ from $w$. Examples of strongly stable sets of words can be obtained by adjoining the empty word to the stable sets of words in the examples   in Section \ref{lumodic}. 
For a  strongly stable set $S$ and a word $w$, set $\Delta_S(w)=  \sum_{w'} w'\otimes w/w'$ where $w'$ runs over
all subwords of $w$ belonging to $S$. It is easy to see that   $\Delta_S$ extends to a coassociative algebra comultiplication in $\WW$. When $S$ is the set of all words, this is the  shuffle comultiplication mentioned in the introduction.

 \section{Hopf algebra  of phrases}\label{oacphras}

The results above suggest  that there may exist  a left-handed   comultiplication in the tensor algebra $T(\VV)$ with leading term $\rho_L$.
  We construct such a comultiplication.

 \subsection{Phrases}\label{modic1}   A {\it phrase of length $k\geq 1$} is a   sequence of $k$   words  in the alphabet $\E$. Some (or all) of these words may be empty. We also allow an {\it empty phrase}  consisting of  0 words and  denoted   $1$.  A more interesting  example of a  phrase: \lq\lq $NIHIL$   $NOVI$ $SUB$ $SOLE$" where   $\E$ is the set of capital Latin  letters. 
  Here words   are separated by  blank spaces and the quotation marks indicate the beginning and the end of the phrase.
 These  conventions, customary in ordinary texts, are  not quite convenient for   mathematical formulas.
 In the formulas we shall indicate the beginning and the end of a phrase by round brackets and   separate words  by   vertical bars. Of course, we assume that the round brackets and the vertical bar are not letters in $\E$. The same example can be re-written as $(NIHIL \V NOVI \V SUB \V SOLE)$.
   By abuse of notation, the   phrase $(w)$ consisting of   one   word $w$ will be also denoted   by $w$.

Let $\MM=\MM(\E)$ be  the free $R$-module freely  generated by the set of phrases in the alphabet  $\E$.   Concatenation of phrases  defines a multiplication in $\MM$ by  
$$(w_1\V w_2 \V... \V w_k) (x_1\V x_2 \V... \V x_{t}) =(w_1\V w_2 \V ... \V w_k\V x_1\V x_2 \V ... \V x_{t}) $$
where $k,t\geq 0$ and $w_1,...,w_k,x_1,...,x_{t}$ are   words.
This 
  makes $\MM$ into an associative algebra with unit $1$ (the empty phrase). This algebra is graded, the grading being   the length of   phrases.
  
  In this section we shall focus on the subalgebra $\PP=\PP(\E)$ of $\MM$ additively generated by phrases formed   by non-empty words. Thus $(w_1\V w_2 \V... \V w_k)\in \PP$ iff $k=0$ or $k\geq 1$ and all the words $w_1,...,w_k$ are non-empty. Note that $1\in \PP$. The algebra $\PP$ is graded: $\PP=\PP^0\oplus \PP^1\oplus ...$ where $\PP^k$ is the module additively generated by phrases of length $k$. The inclusion $\VV=\PP^1\hookrightarrow \PP$ extends to an algebra isomorphism $T(\VV)\to \PP$ and we    identify  the tensor algebra $T(\VV)$ with $ \PP$ along this isomorphism. 

Phrases in the alphabet $\E$ can be viewed as   words in the extended alphabet $\E\amalg \{\V\}$ obtained by adjoining the new letter $\V$ to $\E$. However,  multiplication of phrases is different from multiplication of words
in this extended alphabet due to the  additional symbol $\vert$ between $w_k$ and $x_1$ in the formula above.

\subsection{Comultiplication in $\PP$}\label{pom} Fix a  stable set of words $L$. We define here    a  comultiplication $\Delta=\Delta_L$ in $ \PP$. We begin with preliminary definitions.

Let $w$ be a word   of length $m\geq 1$. A {\it cut} of $w$ of length $k\geq 0$ is a  sequence of $2k$ integers  $1\leq i_1< j_1<  i_2<  j_2< ... < i_k<  j_k \leq m+1$ such that     the $k$ words $w_{i_1,j_1},...,w_{i_k,j_k}$ belong to $L$. For $k=1$, we additionally require that $(i_1,j_1)\neq (1,m+1)$. A cut of length 1 is nothing but a simple cut in terminology of Section \ref{prelilumodic}.  Every word has a unique {\it empty cut} $\emptyset$  of length $0$.  

To indicate that $c=(i_1,j_1,...,i_k,j_k)$ is  a cut of $w$ we   write $c\ll  w$. 
We also write   $\#c= \{1,2,..., k\}$.  For $u\in \#c$, set  $w^c_u=w_{i_u,j_u}$. The factors $w^c_1,...,w^c_k$ of $w$ are called   {\it $c$-factors}.  Finally we define a phrase $$l_c(w)=\prod_{u\in \#c} (w^c_u)= (w^c_{1}\V w^c_{2}\V ...\V w^c_{k})$$ and a     word
$$r_c(w)= w_{1,i_1}  w_{j_1,i_2}w_{j_2,i_3}... w_{j_{k-1},i_k} w_{j_{k},m+1}.$$  
To specify    a cut   of $w$ of length $k$ it is enough to specify    $k$ non-overlapping  proper  factors of $w$ belonging to $L$ and such that   consecutive factors are  separated by at least one letter.    These factors form  the phrase $l_c(w)$. Deleting them     from $w$ we obtain   the word $r_c(w)$.

 Set
$$\Delta(w)=w\otimes 1 + \sum_{c \ll  w} l_c(w) 
\otimes r_c(w)  $$
where $c$ runs over all cuts of $w$.   
Note that the term   corresponding to the empty cut    is $1\otimes w$. The mapping $w\mapsto \Delta(w)$  extends uniquely to an algebra homomorphism $\Delta=\Delta_L:\PP\to \PP^{\otimes 2}$.

 For example, for a word   $w$   of length 1, we have  $\Delta(w)= 
w\otimes 1 +  1\otimes w $. If $L$ is the set of all non-empty words  and    $A, B, C\in \E$,  then 
$$\Delta(AB)= AB \otimes 1+ 1\otimes AB+A\otimes B + B\otimes A,$$ 
$$\Delta(ABC)= ABC \otimes 1+ 1\otimes ABC+A\otimes BC + B\otimes AC+C\otimes AB + AB \otimes C+ BC \otimes A 
+(A\V C)\otimes B  ,$$
$$\Delta(AB\V C)= \Delta (AB) \Delta(C)= (AB\V C)\otimes 1+ C\otimes AB +(A\V C)\otimes B
+(B\V C)\otimes A $$
$$+AB \otimes C +1\otimes (AB\V C)+ A\otimes (B\V C) + B \otimes (A\V C).$$

\begin{theor}\label{comumu} The pair $(\PP,\Delta)$   is  a  left-handed graded bialgebra with leading term   $\rho_L$. \end{theor}
                     \begin{proof} The    left-handedness and the claim  concerning the leading term  follow  directly  from the definitions. The only non-obvious assertion  is the   coassociativity of $\Delta$.
It suffices to prove that $(\id\otimes \Delta)\Delta (w)=(\Delta\otimes \id)\Delta (w)$ for any word $w$. Set
 $$\Theta (w)=\Delta(w)-w\otimes 1= \sum_{c \ll  w} l_c(w) 
\otimes r_c(w). $$
The mapping $w\mapsto \Theta (w)$  defines an $R$-linear homomorphism $\Theta:\VV\to \PP \otimes \VV$. We have 
$$(\Delta \otimes \id) \Delta(w)=(\Delta \otimes \id) \Theta(w)+\Delta(w)\otimes 1 
=(\Delta \otimes \id) \Theta(w)  + \Theta(w)\otimes 1 +w\otimes 1 \otimes 1.$$
Similarly,
$$(\id \otimes \Delta) \Delta(w)=(\id \otimes \Delta) \Theta(w)+ (\id \otimes \Delta)(w\otimes 1 )= 
(\id \otimes \Theta) \Theta(w)+\Theta(w)\otimes 1+w \otimes 1\otimes 1.$$
Comparing these expressions we conclude  that  it is enough to prove that $(\Delta \otimes \id) \Theta(w) =(\id \otimes \Theta) \Theta(w)$.
If follows from the definitions that
$$
(\id \otimes \Theta) \Theta(w)=\sum_{c \ll  w} \,  \sum_{e\ll  r_c(w) }   l_c(w)\otimes l_{e} (r_c(w))\otimes r_{e} (r_c(w)).$$ 
We  compute the right-hand side as follows.

For integers $i<j$, set $[i,j)=\{s\in \ZZ\V  i\leq s<j\}$. For   cuts 
 $d=(i_1,j_1,...,i_k,j_k),d'=(i'_1,j'_1,...,i'_{k'},j'_{k'})$ of $w$,  
   write $d\subset  d'$ if  
  $\cup_{u=1}^k [i_u,j_u) \subset  \cup_{v=1}^{k'} [i'_v,j'_v) $.  Suppose that $ d\subset d'$. We say that the  index $v\in \{1,...,k'\}$ is {\it special}   if    $i'_v=i_u$ and $j'_v=  j_u$ for some $u=1,...,k$.   If $v$ is non-special, then 
   clearly 
   $[i'_v,j'_v) \nsubseteq \cup_{u=1}^k [i_u,j_u)$. We can obtain a cut $\tilde d'$ of $r_d(w)$
by deleting    all   $d$-factors  from   the  $d'$-factors of $w$ numerated by  non-special indices.
That the $\tilde d'$-factors of  $r_d(w)$ belong to $L$ follows from the stability of $L$. 
The cut $\tilde d'$ is empty iff $d=d'$.
Moreover,  the formula $(d,d')\mapsto (d,  \tilde d')$  establishes 
  a bijective correspondence   between   pairs $( d\ll  w,d'  \ll  w)$ with $d\subset d'$  and pairs  $( c \ll  w, e\ll  r_c(w) )$. Therefore
  $$
(\id \otimes \Theta) \Theta(w)=\sum_{d\ll  w ,  d'\ll  w,  d\subset  d' }  
 l_d(w) \otimes    l_{\tilde d'} (r_d(w))  \otimes  r_{\tilde d'} (r_d(w))=\sum_{d\ll  w ,  d'\ll  w,  d\subset  d' }  
 l_d(w) \otimes    l_{\tilde d'} (r_d(w))  \otimes  r_{d'}(w) . $$
We claim that the right-hand side is   equal to $(\Delta \otimes \id) \Theta(w)$. (The proof of this  does  not use the stability of $L$). 
If follows from the definitions that
$$
(\Delta \otimes \id) \Theta(w)=\sum_{c \ll  w}  \prod_{u\in \#c}   \Big (w_u^c\otimes 1  + \sum_{e_u\ll  w^c_u }   l_{e_u} (w^c_u) \otimes r_{e_u} (w^c_u)\Big )
\otimes r_c (w)$$  
$$=\sum_{c \ll  w}   \sum_{I\subset \#c} \,\, \sum_{\{e_u\ll   w^c_u \}_{u\in \#c-I}}   \Big ( \prod_{v\in I } w^c_v  \prod_{  u\in \#c-I}
   l_{e_u} (w^c_u)\Big ) \otimes
\prod_{  u\in \#c-I} r_{e_u} (w^c_u) 
\otimes r_c(w).$$ Here  all the products are   ordered  in accordance with the natural order   in $\#c$. For example, if $\#c=\{1,2,3\}$ and $I=\{1,3\}$, then the term in the big round brackets on the right-hand side is  $ w^c_1\, l_{e_2} (w^c_2) \, w^c_3$. 

Given a cut $c=(i_1,j_1,..., i_k,j_k)$  of $w$ and $u\in \#c$,  every cut $e_u$ of   $w^c_u$ yields a cut $\hat e_u$ of $w$   by adding $i_u-1$ to all terms of $e_u$. A set $I\subset \#c$ gives rise to a cut $\hat I$ of $w$ formed by the indices $\{i_v, j_v\}_{v\in I}$. 
With a 
  tuple  $(c\ll  w, I\subset \# c, \{e_u\ll  w^c_u\}_{u\in \#c-I})$ we   associate 
a pair $(d,d')$ where $d'=c$ and $d$ is the cut   of $w$ obtained as the union of $\hat  I$ with  all $\{\hat e_u   \}_{u\in \#c-I} $.   This  establishes  
  a bijective correspondence between   such tuples  
and   the pairs $(d   \ll  w, d'  \ll  w)$ with $ d\subset  d'$. The corresponding  terms  in the expansions of $(\id \otimes \Theta) \Theta(w)$ and $(\Delta \otimes \id) \Theta(w)$ are equal.  
    Therefore $ (\id \otimes \Theta) \Theta(w)=(\Delta \otimes \id) \Theta(w)$.   \end{proof}

\begin{corol}\label{cor11}   The bialgebra  $(\PP ,\Delta)$ is    a Hopf algebra.
\end{corol}
                     \begin{proof} The augmentation $ \varepsilon: \PP \to R$ sending all non-empty phrases to 0 and sending 1 to 1  is a counit of  $\Delta$. It remains to show the existence of  an antipode. Consider the unique $R$-linear endomorphism $s$ of $ \PP$ such that $s (1)=1$, $s(ab)=s(b) s(a)$ for all $a,b\in \PP$  and  the value  of $s$ on  words  is defined  by induction on  the length as follows:  for a word $w$ of length 1, set $s(w)=-w$; for a word $w$ of    length $\geq 2$, set
$$s(w)=-w-\sum_{c\ll  w, c\neq \emptyset}     l_c(w)   \, s(r_c(w)) \in \PP$$
where we use that $  r_c(w) $ is shorter than $w$.    These formulas imply  that $\mu (\id \otimes s)\Delta (w) =\varepsilon(w)$ where $\mu$ is multiplication in $\PP$. In other words,  
$s$ is a right  inverse of $\id :\PP\to \PP$ with respect to the (associative) convolution product $\bullet$ in $\End_R(\PP)$ defined by $f\bullet g= \mu (f\otimes g)\Delta$ for $f,g\in \End_R(\PP)$.   Similar inductive formulas show that $\id $ has a left inverse
$s'\in \End_R(\PP) $ and then $s'=s' \bullet (\id\bullet s) =(s' \bullet  \id)\bullet s=s$. Therefore $s$ is   an antipode for  $(\PP, \Delta)$. \end{proof}
 
 \subsection{Phrase indicators}\label{dhaa}  By   {\it phrase indicators} we mean $R$-linear homomorphisms $\PP\to R$.  
By the general theory of bialgebras, the  comultiplication $\Delta_L$ in $\PP$ induces an associative multiplication $\circ_L$   in  the module of phrase indicators $\PP^*=\Hom_R(\PP,R)$.   The product $f\circ_L  g$  may distinguish phrases indistinguishable by $f,g\in \PP^*$. For example, let $L$ consist  of all non-empty words and let  $\ell$ and $f_B$ be the phrase indicators counting the number of   words in a phrase and the number of appearances of   $B\in \E$ in a phrase, respectively. Then the values of $\ell \circ_L f_B$ on the 1-word phrases $(ABC)$ and $(ACB)$ (where $A,B,C$ are distinct letters in $\E$) are $4$ and $3$, respectively.

 The  additive homomorphism $\Theta:\VV\to \PP \otimes \VV$ constructed 
in the    proof of Theorem \ref{comumu}   makes $\VV$ into a  comodule over the bialgebra $(\PP, \Delta)$. This induces a right action of the algebra $\PP^*$ on $\VV$.  Using the antiautomorphism of $\PP^*$ induced by the antipode in $\PP$, we can transform the right  action of $\PP^*$ into a left action. 
The    right and left  actions of a phrase indicator $f\in \PP^*$ on $\VV$ depend  only on the values of $f$ on phrases with all words in $L$. This can be formalized as follows.  Let $\PP_L$ be the subalgebra of $\PP$ additively generated by the phrases whose all words belong to $L$ (including the empty phrase 1). 
 It is clear   that $\Delta_L(\PP_L)\subset \PP_L \otimes \PP_L $  so that  $\PP_L$ is a Hopf subalgebra of $\PP$. Clearly,  $\Theta(\VV)\subset \PP_L\otimes \VV$. In this way $\VV$ acquires the structure of  a   comodule over  $\PP_L$. The   actions of the algebras $\PP^*, (\PP_L)^*$ on $\PP$ are compatible    via the algebra homomorphism  $\PP^*\to (\PP_L)^*$ induced by the inclusion $\PP_L\subset \PP$.

\subsection{Dual Hopf algebra}\label{pppldf}  We can define a Hopf algebra dual to $\PP$.  Consider the algebra 
$\PP^*$ with multiplication induced by $\Delta_L$ and quasi-comultiplication  induced by multiplication in $\PP$.
 Consider the   embedding $\delta: \PP\hookrightarrow \PP^*$  mapping a phrase $p$ into the phrase indicator $\delta_p$ whose value on $p$ is 1 and whose value on all other phrases is $0$.  It is easy to see that $\delta(\PP)$ is a subalgebra of $\PP^*$. 
 In this way the module $\PP$ acquires a new associative multiplication   $\circ_L$. The       quasi-comultiplication  in $\PP^*$   induces  a genuine comultiplication in     $\PP$ transforming a phrase
 $(w_1\V...\V w_k)$ into $\sum_{i=0}^k (w_1\V...\V w_i) \otimes (w_{i+1}\V...\V w_k)$. This makes the algebra $(\PP, \circ_L)$ into a   Hopf algebra. By its very definition, it is dual to $(\PP, \Delta_L)$.

 For completeness, we  describe multiplication  $\circ_L$ in $\PP$   explicitly. 
 For a   phrase $p=(w_1\V ...\V w_k)$ and a non-empty  word $y$,  set $p\ast y=0$ if at least one of the words $w_1,...,w_k$ does not belong to $L$. If $w_1,...,w_k\in L$, set
 $$p\ast y=\sum_{ y=x_1\cdots x_{k+1}} x_1w_1x_2w_2... w_k x_{k+1}\in \VV\subset \PP$$
where the sum runs over all  sequences of words   $x_1,  ...,x_{k+1}$ (some of them possibly empty)  such that $y=x_1 \cdots x_{k+1}$. In particular if $p=1$   (that is if $k=0$), then $p\ast y=y $. For   empty word $y$, set     $p\ast y=0$ if $k\neq 1$ and $p\ast y=w_1$ if $k=1$.  
  
Given   a non-empty   phrase $p=(w_1\V ...\V w_k)$, denote by $S (p)$ the set of all finite sequences of phrases 
   $p_1,..., p_t$ (some of them possibly empty) such that $p=p_1p_2...p_t$.   
Denote by $ W (p)$ the set of all  finite sequences of words     obtained by inserting    empty words in the sequence $w_1,w_2,..., w_k$.

It is easy to verify that $1 \circ_L p=p\circ_L 1=p$ for all $p\in \PP$. For     non-empty phrases $p, q$,
 $$p\circ_L q= \delta^{-1} (\delta_p \delta_q)=\sum_{t\geq 1}\,\, \sum_{ (p_1,...,p_t)\in S (p),  (y_1,...,y_t)\in  W (q)}\, 
  \prod_{i=1}^t p_i \ast y_i \in \PP.$$
The right-hand side contains only a finite number of non-zero terms 
  because $p_i\ast y_i=0$ unless $p_i$ and/or $y_i$ are non-empty.

 \subsection{Example}\label{derexidifr} Let the alphabet $\E$ consist  of   one letter $A$.   Then $\PP$ is a free associative (non-commutative)  unital algebra over $R$ freely generated by the words $A, A^2=AA, A^3=AAA,...$ (caution: $A^2$ is not the square of $A$ in $\PP$).  The comultiplication $\Delta_L$ in $\PP$ corresponding  to $L=\{A^m\}_{m\geq 1}$ is computed by
 $$\Delta_L(A^m)=A^m\otimes 1 + 1\otimes A^m+ \sum_{m>m_1\geq 1} (m+1-m_1) \,A^{m_1} \otimes A^{m-m_1}$$
 $$ +\sum_{k\geq 2}\, \sum_{m_1,...,m_k\geq 1, 
  m_1+...+m_k \leq m-k+1} { {m+1-m_1-...-m_k}\choose  { k}} \, A^{m_1}   \cdots A^{m_k} \otimes A^{m-m_1-...-m_k}.$$
 If $L$ consists of the  words  whose length is divisible by a given   integer $N\geq 1$, then the formula for $\Delta_L$ is the same with the  restriction that $m_1,...,m_k$ are divisible by $N$. 
 
 \subsection{Functoriality}\label{funcpppldf}  Any mapping $\alpha$ from an alphabet $\E$ to an alphabet $\E'$
 extends to words  letter-wise. Denote the resulting mapping  by   $\tilde \alpha$. A  stable set of words $L'$ in the alphabet $\E'$ gives rise to a stable set of words $L=\tilde \alpha^{-1}(L')$ in the alphabet $\E$. It is clear that $\tilde \alpha$ induces  an   $R$-homomorphism   $  (\VV(\E), \rho_{L})\to (\VV(\E'),\rho_{L'})$  
 of pre-Lie coalgebras. The latter extends by multiplicativity to a homomorphism $ (\PP(\E), \Delta_{L})\to (\PP(\E'),\Delta_{L'})$ 
 of Hopf algebras. For example, if $L$ is the set of all non-empty words in the alphabet $\E$, then any permutation of $\E$ induces an automorphism of the pre-Lie coalgebra $(\VV(\E), \rho_{L})$ and an automorphism of the Hopf  algebra $(\PP(\E), \Delta_{L})$.

  \section{Coalgebra   of words: second  construction}\label{seccoac}

Fix from now on a  mapping $\mu : \E\times \E \to R$.  We derive from $\mu$ a   pre-Lie comultiplication $\rho=\rho_\mu$ in   the $R$-module $\WW=\VV\oplus R\phi$.

 \subsection{Simple inscriptions}\label{fginscdom}    Let $w$ be a word   of length $m\geq 0$. A  {\it simple   inscription}   in $w$ is a  pair $i,j \in  \{1,2,...,m\}$ with $i<j$. To indicate that $a =(i,j)$ is  a  simple inscription in   $w$ we   write $a\dag  w$.         
  Consider the word  $ l_a(w)= w_{i+1,j}$ of length $j-i-1$ and the word $ r_a(w)=w_{1,i}   w_{j+1,m+1} $  of length $  m+i-j-1$.   Set $\mu(w\vert_a)= \mu (w (i ), w(j)) \in R$.     For example, if $w=ABABB$ and $a=\{1, 4 \}$, then   $r_a(w)=B$,  
$l_a(w)=  BA$, and  $\mu(w\vert_a)= \mu (A,B)$. 

 Set
 $$\rho(w)=  \sum_{a\dag w} \mu (w\vert_a) \, l_a(w)\otimes r_a(w) \in \WW \otimes \WW  $$ where $a$ runs over all simple inscriptions in  $w$.
If $w $ is an empty word or  a 1-letter word, then $w$ has no simple inscriptions so that  $\rho(w)=0$. 
    Extending $\rho$ by linearity, we obtain a comultiplication $\rho:\WW \to \WW \otimes \WW $.

\begin{theor}\label{domumu}     $\rho$   is a  pre-Lie comultiplication in $\WW$.
\end{theor}
                     \begin{proof}  Consider a   word $w$ of length $m $ and  two simple   inscriptions $a=(i,j)$ and $b=(i', j')$    in $w$.  We write $a< b$ if $j<   i'$. In this case set
		     $$w(a,b)= \mu (w\vert_a) \, \mu (w\vert_b) \, l_{b}(w) \otimes l_a(w)\otimes  w_{1,i}  w_{j+1,i' } w_{j'+1,m+1}   \in \WW^{\otimes 3}.$$
We write $a>b$ if $i>    j'$ and set then  
		     $$w(a,b)= \mu (w\vert_a) \, \mu (w\vert_b)\, l_{b}(w) \otimes l_a(w)\otimes  w_{1,i'}  w_{j'+1,i }  w_{j+1,m+1}   \in \WW^{\otimes 3}.$$	Note that $a>b$ iff $b<a$ and then $w(a,b)= P^{1,2}  (w(b,a))$.	
	 We write      
  $a\sqsubset b$ if   $i < i'< j' < j$.     Set then  $$w(a,b)= \mu (w\vert_a) \, \mu (w\vert_b)\, l_{b}(w) \otimes   w_{i +1,i'} w_{j'+1,j}  \otimes  r_a(w)  \in \WW^{\otimes 3}.$$

 We expand $(\id\otimes \rho )\rho (w)$ from definitions:
 $$ (\id\otimes \rho )\rho (w)= \sum_{b \dag w}	  \mu (w\vert_b) \,l_b(w)\otimes \rho  (r_b(w))$$
 $$
=\sum_{ b \dag w }\,  \mu (w\vert_b)\,	\sum_{ a \dag r_b(w) } \mu (r_b(w)\vert_a)  \,  l_b(w)\otimes l_{a} ( r_b(w))\otimes r_{a} ( r_b(w)).$$
  For  $b=(i', j')\dag w$, we can describe all   simple inscriptions in  $r_b(w)=w_{1,i'}   w_{j'+1,m+1}$  as follows.
A simple inscription   $a\dag w$ such that    $a< b$   is automatically a  simple inscription in $r_b(w)$. A simple inscription $a=(i,j)$ in $ w$ such that     $a> b$
  yields a  simple inscription   in $r_b(w)$  by subtracting $j'+1-i'$ from both $i$ and $j$. A simple inscription
$a=(i,j)$ in $w$ such that $a \sqsubset b$ yields a  simple inscription $(i, j-(j'+1-i'))$ in $r_b(w)$.  It is clear that every simple inscription in $r_b(w)$    arises  in exactly one of these 3 ways from a certain   $a\dag w$.  The corresponding term in the expansion of  $(\id\otimes \rho )\rho (w)$ is $w(a,b)$.
Therefore
 $(\id\otimes \rho )\rho (w)=x+y+z$ where 
 $$x=\sum_{a,b\dag w, a<b} w(a,b),\,\,\, \,\,\,y=\sum_{a,b\dag w, a>b} w(a,b),\,\,\,\,\,\, z=  \sum_{a,b\dag w,  a \sqsubset b} w(a,b).$$ 
A similar (in fact easier) computation shows that  $(\rho \otimes \id) \rho (w)=z$. 
Therefore
$$\tilde \rho (w)=((\id\otimes \rho )\rho -(\rho \otimes \id) \rho ) (w)= x+y .$$
By the remarks above,  $y=   P^{1,2}  (x)$ so that $\tilde \rho (w)=x+y$ is   invariant under $P^{1,2} $. \end{proof}

 \subsection{Extended word indicators}\label{wordomph} By the general theory exposed in Section 2, 
 the pre-Lie comultiplication $\rho_\mu$ induces a pre-Lie multiplication $\star_\mu$  and a Lie bracket $[f,g]_\mu=f\star_\mu g- g\star_\mu f$ in the module $\WW^*=\Hom_R(\WW,R)$. The elements of $\WW^*$ are called  {\it  extended  word indicators}.  
 The module $\WW^*$ admits a decreasing filtration $\WW^*=   \WW^{*(0)}\supset \WW^{*(1)}\supset ...$ where $\WW^{*(m)}$  consists of the   indicators annihilating all words of length  $\leq  m-1$.  It is clear that $\WW^{*(m)} \star_\mu \WW^{*(n)}\subset \WW^{*(m+n+2)}$  and   $[\WW^{*(m)} , \WW^{*(n)}]_\mu\subset \WW^{*({m+n+2})}$ for all $m,n $. 
 This   implies that   $\WW^*= \projlim_m 
(\WW^*/\WW^{*(m)})$
is a  projective limit of nilpotent Lie algebras.

 Recall the   Lie algebra action of $\WW^*$   on $\WW$ induced by $\rho_\mu$.    For   $f\in \WW^*$ and a   word $w$,
 $$fw=-\sum_{a\dag w} \,  \langle l_a(w), f\rangle \,\mu (w\vert_a)\,r_a(w) \in \WW.$$
  Consider the filtration $0=\WW_{-1}\subset R\emptyset =\WW_0\subset \WW_1\subset   ...$ of $\WW$
 where $\WW_m$ is generated by the words of length $\leq m$. It is clear that   $f\WW_m\subset \WW_{m-2}$ for all $m$ and all  $f\in \WW^*$. This implies that   the action of $\WW^*$ on $\WW_m$ is nilpotent for all $m$.

Consider the  embedding $\delta: \WW\hookrightarrow \WW^*$  mapping a word $w$ into the extended word indicator $\delta_w$ whose value on $w$ is 1 and whose value on all other words is $0$. If $\E$ is finite or more generally if $\mu$ takes non-zero values only on a finite subset of $\E\times \E$, then  the image of $\delta$ is closed under $\star_\mu$. This induces a pre-Lie  multiplication   on   $\WW$. We  leave it to the reader to give an explicit formula for it.

If $R\supset \QQ$, then  the Campbell-Baker-Hausdorff formula defines   a
group
multiplication in $\WW^*$ as in Section \ref{succ3}.  
The resulting pro-nilpotent group
  $\Exp \WW^*$  is
the \lq\lq Lie group" with Lie algebra $\WW^*$.  The  Lie algebra action of $\WW^*$ on  $\WW$ induced by $\rho$ integrates  to a group action of $\Exp \WW^*$ on  $\WW$ as in Section \ref{succ3}.

\subsection{Example}\label{exann} 
  Let $\mu :\E\times \E\to R$ send  a pair $(A,B)$ to $1$ if $A=B$ and to 0 if $A\neq B$.   Let   $w=ABACBA$ with $A,B,C\in \E$. Then
$$\rho_\mu (w)= B\otimes CBA + CB \otimes AB+ BACB \otimes \phi
+   AC\otimes  AA .$$
For any word indicator $f$,
$$fw= -\langle  B  , f\rangle\, CBA -\langle  CB  , f\rangle\, AB
 -\langle  BACB  , f\rangle\, \phi  -  \langle AC  , f\rangle\, AA.$$
For example, let  the indicator  $ f_A$ compute  the total number of occurencies of the letter $A$ in a word. Then
  $f_A w=-\phi- AA$, $f_Bw=-CBA- AB- 2\phi $, and  
$f_C w=-  AB -\phi  -  AA$.  If $\ell$ is the indicator   computing  the length of a word, then 
$\ell w=-CBA-2 AB -4\phi   -   2 AA$.
 
 \section{Hopf algebra  of phrases: second construction}\label{secoacphras}

 \subsection{Comultiplication   $\Delta_\mu$}\label{commtepom}  Recall the algebra of phrases $\MM$ defined in  Section \ref{modic1}. The   inclusion $\WW\hookrightarrow \MM$ as 1-word phrases extends to an isomorphism  of the tensor algebra $T(\WW)$ onto $ \MM$. We shall identify   $T(\WW)$ with   $\MM$.  The results above suggest  that there may exist  a left-handed   comultiplication in  $\MM $    with leading term $\rho_\mu$.
 We define such  a  comultiplication $\Delta=\Delta_\mu$ in $ \MM$. 

  Let $w$ be a word   of length $m\geq 0$.  By an {\it inscription} in $w$ we shall mean a subword of $w$ of even length.  More precisely, an       inscription  in $w$  of length $k\geq 0$ is a  set   $\alpha\subset \{1,2,...,m\}$ consisting of $2k$ elements. We shall list  these elements in the increasing order and write $\alpha =(i_1,j_1,...,i_k,j_k)$ where $i_1<j_1<...<i_k<j_k$.    Every word has a unique {\it empty inscription} of length  $ 0$. To indicate that $\alpha =(i_1,j_1,...,i_k,j_k)$ is  an inscription in   $w$ we   write $\alpha\ddag  w$.
 Set $\#\alpha=\{1,2,...,k\}$ and $supp (\alpha)=\cup_{u\in \#\alpha} [i_u, j_u]$ where $[i,j]=\{k\in \ZZ\,\vert\, i\leq k\leq j\}$.  With each $u\in \#\alpha$, we associate the word $w^\alpha_u=w_{i_u+1,j_u}$. It is empty iff $i_u+1=j_u$.
We define a  phrase  
  $$l_\alpha(w)=\prod_{u\in \#\alpha} w^\alpha_u= (w^\alpha_1 \V w^\alpha_2\V... \V w^\alpha_k)\in \MM .$$
  Clearly $l_\alpha(w)=1$  iff $\alpha$ is void. We also associate with $\alpha$
 an element    of the ground ring 
 $$\mu (w, \alpha)=\prod_{u\in \#\alpha}  \mu(w(i_u), w(j_u) )  $$
 and a word
 $$r_\alpha(w)= w_{1, i_{1}} w_{j_{1}+1, i_{ 2}} \cdots  w_{j_{k-1}+1, i_{k}} w_{j_{k}+1, m+1} . $$
  If $\alpha=\emptyset$, then $l_\alpha(w)=1$, $\mu (w, \alpha)=1$, and  $r_\alpha(w)=w$.
 
Set
 $$ \Delta (w)=w\otimes 1   +\sum_{\alpha\ddag w} \,  \mu(w,\alpha)\,
 l_{\alpha}(w) 
 \otimes r_{\alpha}(w)  $$
   where $\alpha $   runs   over all    inscriptions in  $w$.
Note that  the term   corresponding to    $\alpha =\emptyset$  is $1 \otimes w$.  
 
 The mapping $w\mapsto \Delta(w)$  extends uniquely to an algebra homomorphism $\Delta :\MM\to \MM{\otimes}\MM$.

\begin{theor}\label{coclnewcomumu} The pair $(\MM,\Delta)$   is  a  left-handed graded bialgebra with leading term   $\rho_\mu$. \end{theor}
                     \begin{proof}   The only non-obvious assertion  is the   coassociativity of $\Delta$.
It suffices to prove that $(\id\otimes \Delta)\Delta (w)=(\Delta\otimes \id)\Delta (w)$ for any word $w$. Set
 $$\Theta (w)=\Delta(w)-w\otimes 1=   \sum_{\alpha\ddag w}   \mu(w,\alpha)\, 
 l_{\alpha}(w)
 \otimes  r_{\alpha}(w) . $$
 The mapping $w\mapsto \Theta (w)$  defines an $R$-linear homomorphism $\Theta:\WW\to \PP \otimes \WW$. A computation similar to the one in the proof of Theorem \ref{comumu} shows  that  it is enough to prove that $(\Delta \otimes \id) \Theta(w) =(\id \otimes \Theta) \Theta(w)$. 

 If $w $ has length $0$ or $1$, then  $\Theta(w) = 1\otimes w$ and $(\Delta \otimes \id) \Theta(w) =1\otimes 1 \otimes w=(\id \otimes \Theta) \Theta(w)$.
Suppose from now on that $w$ has length $\geq 2$.   For   inscriptions $\beta $ and $\eta $ in $w$, we write $\beta\lhd \eta$ if    $supp(\beta)\cap \eta=\emptyset$.
Striking out from $w$ all letters numerated by elements of the set $supp(\beta)$ we obtain the word $r_\beta(w)$.  If $\beta\lhd \eta $, then the letters of $w$ numerated by elements of $\eta$ survive in  $r_\beta(w)$ and form  
an inscription  in $r_\beta(w)$  denoted $\eta/\beta$.    The formula $(\beta, \eta)\mapsto (\beta,     \eta/\beta)$  establishes
  a bijective correspondence   between   pairs $( \beta\ddag  w,\eta   \ddag  w)$ such that $\beta\lhd \eta$
     and    pairs  
  $( \alpha \ddag  w, \gamma\ddag  r_\alpha(w) )$.  
Therefore
 $$(\id \otimes \Theta) \Theta(w)=
 \sum_{\alpha \ddag  w} \, \mu (w, \alpha)\,   \sum_{\gamma \ddag  r_\alpha(w) } \, \mu (r_\alpha(w), \gamma)\,  l_\alpha(w)\otimes l_{\gamma} (r_\alpha(w))\otimes r_{\gamma} (r_\alpha(w)) $$
  $$ =\sum_{\beta, \eta \ddag w, \beta\lhd \eta}\,  \mu (w, \beta)\, \mu (w, \eta)\, l_\beta(w)\otimes  l_{  \eta/\beta} (r_\beta(w)) \otimes r_{\eta/\beta} (r_\beta(w)).$$ 
On the other hand, 
 $$
(\Delta \otimes \id) \Theta(w)=   \sum_{\alpha   \ddag  w}  \mu (w, \alpha)\, \Delta (l_{\alpha}(w))\otimes  r_\alpha (w) $$
$$=\sum_{\alpha    \ddag  w}  \mu (w, \alpha)\,  \prod_{u\in \#\alpha}   \Big (w^\alpha_u\otimes 1  +\sum_{\varepsilon_u\ddag w^\alpha_u } \mu( w^\alpha_u, \varepsilon_u)\,\, l_{\varepsilon_u}(w^\alpha_u) \otimes  r_{\varepsilon_u}(w^\alpha_u) \Big ) 
\otimes  r_\alpha (w) $$  
$$= \sum_{\alpha   \ddag  w} \mu (w, \alpha)\, \sum_{I\subset \#\alpha}
 \, \sum_{\{\varepsilon_u\ddag w^\alpha_u \}_{u\in \#\alpha-I}}       \prod_{u\in I } w^\alpha_u  \prod_{  u\in \#\alpha-I}
  \mu( w^\alpha_u, \varepsilon_u)\,\, l_{\varepsilon_u}(w^\alpha_u)   \otimes
 \prod_{  u\in \#\alpha-I} r_{\varepsilon_u}(w^\alpha_u)
\otimes   r_\alpha (w) .$$ Here  the products are   ordered  in accordance with the order of the indices in $\#\alpha$. For example, if $\#\alpha=\{1,2,3\}$ and $I=\{1,3\}$, then the   first tensor factor on the right hand side is  $  \mu( w^{\alpha}_2, \varepsilon_2)\, w^{\alpha}_1\, l_{\varepsilon_2}( w^{\alpha}_2   ) \,  w^{\alpha}_3$. 

Consider an inscription $\alpha=(i_1,j_1,..., i_k, j_k)$ in $w$. Given $u\in \#\alpha$ and an  inscription $\varepsilon_u  $ in   $w^\alpha_u$, we obtain  an inscription  (of the same length) $\hat  \varepsilon_u$ in  $w$
by adding $i_u$ to all terms of $\varepsilon_u  $.   With a 
  tuple  $(\alpha\ddag   w, I\subset \#\alpha, \{\varepsilon_u\ddag w^\alpha_u \}_{u\in \#\alpha-I})$ we   associate 
two inscriptions $ \beta,\eta $ in $w$ by  $\beta=\{i_u, j_u\}_{u\in I}\cup \cup_{u\in \#\alpha-I} \hat  \varepsilon_u  $ and $\eta= \{i_u, j_u\}_{u\in \#\alpha - I}$.   This  defines 
  a bijective correspondence between   such tuples  
and   the pairs $( \beta\ddag  w,\eta  \ddag  w)$ such that $\beta\lhd \eta$. The corresponding  terms  in the expansions for  $(\Delta \otimes \id) \Theta(p)$ and $(\id \otimes \Theta) \Theta(p)$ are equal.  
    Therefore $(\Delta \otimes \id) \Theta(p)=(\id \otimes \Theta) \Theta(p)$.   \end{proof}

\begin{corol}\label{newdfcor11}   The algebra  $\MM $ with comultiplication $\Delta$ is    a Hopf algebra.
\end{corol}
                     \begin{proof} The augmentation $   \MM \to R$ sending all non-empty phrases  to 0 and sending the empty phrase  to 1  is a counit of  $\Delta$. The existence of an antipode is shown as in the proof of Corollary \ref{cor11}.
 \end{proof}

  \subsection{Extended phrase indicators}\label{edrftdhaa}  Homomorphisms $ \MM\to R$ are called {\it extended phrase indicators}.  
 The  comultiplication $\Delta_\mu$ in $\MM$ induces an associative multiplication $\circ_\mu$   in   $\MM^*=\Hom_R(\MM,R)$.   It is easy to give examples showing that   $f\circ_\mu  g$  may distinguish phrases indistinguishable by $f,g\in \MM^*$.

 The  additive homomorphism $\Theta:\WW\to \MM \otimes \WW$ constructed 
in the    proof of Theorem \ref{coclnewcomumu}   makes $\WW$ into a  comodule over the bialgebra $(\MM, \Delta)$. The leading term of $\Theta$ is the pre-Lie comultiplication  $\rho_\mu$ in $\WW$.  The coaction $\Theta$ induces a right action of the algebra of extended phrase indicators $\MM^*$ on $\WW$.  Using the antiautomorphism of $\MM^*$ induced by the antipode in $\MM$, we can transform the right  action of $\MM^*$ into a left action.

\subsection{Dual Hopf algebra}\label{edrpppldf}    If $\E$ is finite or more generally if $\mu$ takes non-zero values only on a finite subset of $\E\times \E$, then  a   Hopf algebra dual to $\MM$ can be constructed as follows. Consider the algebra 
$\MM^*$ with multiplication induced by $\Delta_\mu$ and quasi-comultiplication  induced by multiplication in $\MM$. Consider the   embedding $\delta: \MM\hookrightarrow \MM^*$  mapping a phrase $p$ into the phrase indicator $\delta_p$ whose value on $p$ is 1 and whose value on all other phrases is $0$.   Under our assumptions on $\mu$,  $\delta(\MM)$ is a subalgebra of $\MM^*$. 
 This induces a new associative multiplication   $\circ_\mu$ in $\MM$. The       quasi-comultiplication  in $\MM^*$   induces  a genuine comultiplication in     $\MM$ transforming a phrase
 $(w_1\V...\V w_k)$ into $\sum_{i=0}^k (w_1\V...\V w_i) \otimes (w_{i+1}\V...\V w_k)$. This makes the algebra $(\MM, \circ_\mu)$ into a   Hopf algebra. By its very definition, it is dual to $(\MM, \Delta_\mu)$.
 
 \subsection{Independence of the basis}\label{indefuncpppldf}  Let $M$ be the free $R$-module with basis $\E$.
The mapping $\mu:\E\times \E\to R$ extends to a bilinear form $M\times M \to R$ also denoted $  \mu$. 
 The constructions above produce a pre-Lie comultiplication $\rho_\mu$ in $\WW=T(M)$ and a  Hopf comultiplication $\Delta_\mu$ in $\MM=T(\WW)$. An   inspection of these constructions shows that they are entirely determined by the   form $  \mu:M\times M \to R$ and do not depend on the   basis $\E$ in $M$. This version of the constructions  applies to an arbitrary (not necessarily free) $R$-module $M$  endowed with a bilinear form $\mu: M\times M \to R$. It yields   a pre-Lie coalgebra $(T(M), \rho_{\mu})$ and a  Hopf  algebra $(T(T(M)), \Delta_{\mu})$.

 A  homomorphism of $R$-modules $\psi:M\to M'$ compatible with bilinear forms $  \mu:M\times M \to R,  \mu':M'\times M' \to R$ (so that $\mu=\mu' (\psi\times \psi)$) induces    a homomorphism
 of pre-Lie coalgebras $  (T(M), \rho_{\mu})\to (T(M'),\rho_{\mu'})$   and  a homomorphism
 of Hopf algebras $ (T(T(M)), \Delta_{\mu})\to (T(T(M')), \Delta_{\mu'})$. In particular,   $\mu$-preserving automorphisms of $M$   induce  automorphisms of   $(T(M), \rho_{\mu})$ and   of   $(T(T(M)), \Delta_{\mu})$.

 \subsection{Relations with algebras of trees}\label{relsecoacphras}
  Let $\mu:\E\times \E\to R$ be the mapping sending a pair $(A,B)\in \E\times \E$ to $1$ if $A=B$ and to $0$ if $A\neq B$. We relate the    Hopf algebra
$(\MM, \Delta_\mu)$  with the Hopf algebra of planar rooted trees due to Connes-Kreimer \cite{ck} and Foissy \cite{fo}.
  
  We say that a word $w$ in the alphabet $\E$   is {\it unlaced} if (i) every letter of $\E$ either does not appear in $w$ or appears in $w$ twice and (ii) for any distinct letters $A,B $ appearing in $w$ the word $ABAB$ is not a subword of $w$. The latter condition may be reformulated by saying that $w$  must have   the   form   
 $\cdots A \cdots A \cdots B \cdots B \cdots$ or $\cdots A \cdots B \cdots B \cdots A \cdots$
 or $\cdots B \cdots B \cdots A \cdots A \cdots$.  Let $U\subset \WW$ be the free $R$-module generated by the unlaced words (including the empty word).  It follows from the definitions, that    $\rho_\mu(U)\subset U\otimes U$.  Similarly,  the tensor algebra $ T(U)$ is a Hopf subalgebra of $(\MM=T(\WW), \Delta_\mu)$.
 Restricting  $\Delta_\mu$ to $T(U)$ we obtain a Hopf algebra 
 $(T(U), \Delta_\mu)$. 
 
 Unlaced words  can be described in terms of decorated  planar rooted (finite) trees as follows.
 We say that a tree is {\it decorated} if every its edge is labeled by a letter of  the alphabet $\E$ so that different edges are labeled by different letters.
 Each  decorated  planar rooted tree $\tau\subset \RR^2$ gives rise to a word $w(\tau)$  as follows. 
 Consider a  narrow neighborhood $V$ of $\tau$ in $\RR^2$. If $\tau$ has $n$ edges then the circle $\partial V$ consists of $2n$ arcs
 going closely to edges of $\tau$. Starting at a point   $x\in \partial V$ near  the root of $\tau$ and  moving along $\partial V$ counterclockwise we write down the labels of the corresponding edges of $\tau$
 until the first return to $x$. This gives  $w(\tau)$. It is clear from the definitions that this word is unlaced. For example, if $\tau$ is a point, then $w(\tau)=\phi$. If $\tau$ consists of a single edge labeled with $A$, then $w(\tau)=AA$. If $\tau$ is a $Y$-shaped tree with   3 edges, then $w(\tau)=ABBCCA$ where $A$ is the label of the edge incident to the root and $B,C$ are the labels of the two other edges.
 
 An   induction on the number of edges shows that the formula $\tau\mapsto w(\tau)$ establishes a bijective correspondence between
 decorated planar rooted trees (considered up to ambient isotopy in the plane) and unlaced words. In this way the tensor algebra $T(U)$ can be identified with  the free associative (non-commutative) algebra ${\mathcal{H}}$ generated by decorated planar rooted trees. A   comparison of definitions shows  that under this identification  
  the comultiplication $\Delta_\mu$ in $T(U)$
 coincides with the Connes-Kreimer-Foissy comultiplication in ${\mathcal{H}}$. Note that Connes and Kreimer considered a commutative algebra  generated by    rooted trees (without planar structure). A non-commutative version of their definition  was given by Foissy for planar rooted trees. In his paper,  Foissy decorates   vertices of   trees rather than   edges. However his definitions directly extend to  trees with decorated edges.


\begin{thebibliography}{CJKLS}

 \bibitem[CK]{ck} A. Connes, D. Kreimer, \emph{Hopf algebras, renormalization and noncommutative 
 geometry\/},
 Comm. Math. Phys.  199  (1998),  no. 1, 203--242. 


 \bibitem[Fo]{fo} L. Foissy,  \emph{Les alg\`ebres de Hopf des
  arbres enracin\'es d\'ecor\'es. I\/}.    Bull. Sci. Math.  126  (2002), 
   no. 3, 193--239.


\bibitem[Ge]{ge} M. Gerstenhaber,
 \emph{The cohomology structure of an associative ring\/},
Ann. of Math. (2) 78 (1963),  267--288.

 


 

 \bibitem[Tu]{tu1} V. Turaev,
\emph{Skein quantization of Poisson algebras of loops
on surfaces\/},
Ann. Sci. \'Ecole Norm. Sup. (4) 24 (1991), no. 6, 635--704.

\bibitem[Vi]{vi}  E.B. Vinberg,
\emph{The theory of convex homogeneous cones\/},
Trans. Mosc. Math. Soc. 12 (1963), 340-403; translation from Tr. 
Mosk. Mat. Ob-va 12 (1963), 303-358.


                     \end{thebibliography}
                     \end{document}